\documentclass[12pt]{amsart}
\usepackage{amssymb}
\usepackage{pict2e}
\usepackage{amscd}
\def\Z{\mathbb{Z}}
\def\Q{\mathbb{Q}}
\def\R{\mathbb{R}}

\def\ybar{\bar y}
\def\zbar{\bar z}
\def\Xbar{\bar X}

\def\Gbar{\bar G}
\def\a{\alpha}
\def\iso{\cong}
\def\power#1{{\uparrow}\bigl\{#1\bigr\}}
\def\Bigpower#1{{\uparrow}\Bigl\{#1\Bigr\}}
\def\semidirect{\rtimes}
\def\notdiv{\nmid}
\def\-{\phantom{-}}

\newtheorem{theorem}{Theorem}
\newtheorem{lemma}[theorem]{Lemma}
\newtheorem{corollary}[theorem]{Corollary}
\newtheorem*{question}{Question}
\theoremstyle{remark}
\newtheorem*{remark}{Remark}

\begin{document}
\title{Triangles of Baumslag-Solitar Groups}
\author{Daniel Allcock}
\address{Department of Mathematics\\University of Texas, Austin}
\email{allcock@math.utexas.edu}
\urladdr{http://www.math.utexas.edu/\textasciitilde allcock}
\thanks{Partly supported by NSF grant DMS-0600112.}
\date{July 22, 2010}

\begin{abstract}
Our main result is that many triangles of Baumslag-Solitar groups
collapse to finite groups, generalizing a famous example of Hirsch and
other examples due to several authors.  A triangle of Baumslag-Solitar
groups means a group with three generators, cyclically ordered, with
each generator conjugating some power of the previous one to another
power.  There are six parameters, occurring in pairs, and we show that
the triangle fails to be developable whenever one of the parameters
divides its partner, except for a few special cases.  Furthermore,
under fairly general conditions, the group turns out to be finite and
solvable of class${}\leq3$.  We obtain a lot of information about
finite quotients, even when we cannot determine developability.
\end{abstract}

\maketitle

We study groups $G$ of the form
\begin{equation}
\label{eq-deg-of-G}
G(a,b;c,d;e,f):=
\bigl\langle x,y,z \bigm|
(x^a)^y=x^b, (y^c)^z=y^d, (z^e)^x=z^f
\bigr\rangle,
\end{equation}
where $a,\dots,f$ are nonzero integers.  We prove that $G$ collapses
to a finite solvable group under a mild divisibility condition on the parameters.  The
motivation is that $G$ is a triangle of groups in the language of
\cite{Haefliger} or \cite{Stallings}, with the vertex groups being
Baumslag-Solitar groups.  Polygons of groups are an important means of
constructing groups in geometric group theory; see e.g., \cite{BH},
\cite{Farb-Hruska-Thomas}, and \cite{Wise}.  And the Baumslag-Solitar
groups are famous for their ``pathological'' properties, like being
non-Hopfian and (therefore) non-residually-finite and non-linear.

These groups allow a simple construction (probably the first one) of a
non-developable triangle of groups, because $G(1,2;1,2;1,2)$ turns out
to be trivial.  This is a result of K. Hirsch, reported by
Higman \cite{Higman} and motivated by 
Higman's use of a square of $BS(1,2)$'s to construct a
finitely presented infinite group with no finite quotients.
See also \cite[\S23]{Neumann}.  The observation that it can be regarded as a
non-developable triangle of groups seems to be due to K.~Brown.  Here
the vertex groups are copies of $BS(1,2)$, which is an atypical
Baumslag-Solitar group, since it is solvable.  It is natural to ask
what is really causing the collapse; this led to our more general
non-developability criterion:

\begin{theorem}
\label{thm-non-devlopability}
Regard $a$ and $b$ as partners, and similarly for $c$ and $d$ and
for $e$ and $f$, and suppose one of $a,\dots,f$ divides its partner.  Then
the triangle of groups $G(a,b;c,d;e,f)$:
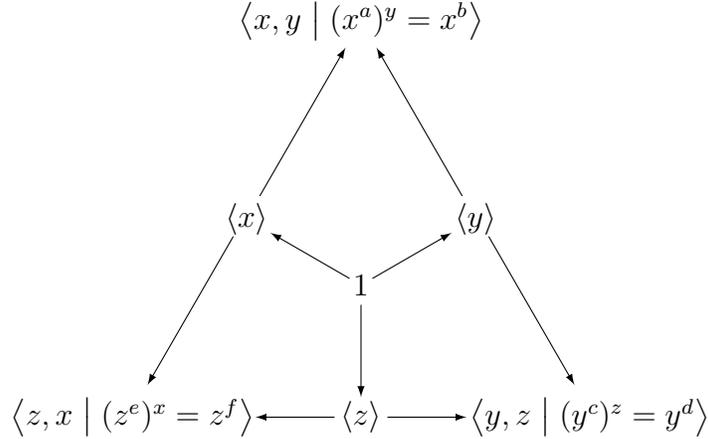
\begin{center}
\unitlength=.5pt
\begin{picture}(200,360)(-100,-150)
\put(0,0){\makebox(0,0)[c]{$1$}}
\put(0,-100){\makebox(0,0)[c]{$\langle z\rangle$}}
\put(86.6,50){\makebox(0,0)[c]{$\langle y\rangle$}}
\put(-86.6,50){\makebox(0,0)[c]{$\langle x\rangle$}}
\put(0,200){\makebox(0,0)[c]{$\bigl\langle x,y\bigm|(x^a)^y=x^b\bigr\rangle$}}
\put(173,-100){\makebox(0,0)[c]{$\bigl\langle y,z\bigm|(y^c)^z=y^d\bigr\rangle$}}
\put(-173,-100){\makebox(0,0)[c]{$\bigl\langle z,x\bigm|(z^e)^x=z^f\bigr\rangle$}}
\put(0,-15){\vector(0,-1){70}}
\put(8.66,5){\vector(866,500){60}}
\put(-8.66,5){\vector(-866,500){60}}
\put(20,-100){\vector(1,0){60}}
\put(-20,-100){\vector(-1,0){60}}
\put(96.6,36.7){\vector(500,-866){65}}
\put(-96.6,36.7){\vector(-500,-866){65}}
\put(76.6,67.3){\vector(-500,866){65}}
\put(-76.6,67.3){\vector(500,866){65}}
\end{picture}
\end{center}
is not developable, except in the special cases
\begin{align}
\label{eq-special-case-1}
&G(a,-a;c,-c;e,-e),\\
\label{eq-special-case-2}
&G(a,\-b;c,\-c;e,\-e),\\
\label{eq-special-case-3}
&G(a,\-b;c,\-c;e,-e)\hbox{, $a\equiv b\mod2$},\\
\label{eq-special-case-4}
&G(a,\-b;c,-c;e,\-e)\hbox{, $e$ even},\\
\label{eq-special-case-5}
&G(a,\-b;c,-c;e,-e)\hbox{, $e$ even and $a\equiv b\mod2$},
\end{align}
all of which are developable.
\end{theorem}
We remind the reader that a triangle of groups is called developable
if each of its vertex groups injects into the direct limit of the
diagram, which in this case is $G(a,b;c,d;e,f)$.  We will be informal
and say that the group is developable when we mean that the triangle
is.  In the list of special cases we have left implicit other cases
obtained from these by ``trivial'' transformations.  These are cyclic
permutation of the three pairs (corresponding to cyclic permutation of
$x,y,z$), exchange of one of $a,\dots,f$ with its partner
(corresponding to inverting one of $x,y,z$), and simultaneous negation
of one of $a,\dots,f$ and its partner (corresponding to inverting a
relation).  We will apply these ``moves'' freely when it is
convenient.

Of course,
theorem~\ref{thm-non-devlopability} begs the question:

\begin{question}
If none of $a,\dots,f$ divides its partner, is $G(a,b;c,d;e,f)$ ever
developable?  always developable?
\end{question}

Our work generalizes results of Post \cite{Post}, who showed finiteness
when $e=1$ and the other parameters satisfy mild inequalities.  His
paper followed work by Mennicke \cite{Mennicke} and
Wamsley \cite{Wamsley} concerning the case $a=c=e=1$; see also Johnson and
Robertson \cite{JR} and most recently Jabara \cite{Jabara}.  The main claim of
Neumann \cite{Neumann-broken} is that $G$ is infinite if $2\leq a\leq|b|$, $2\leq
c\leq|d|$ and $2\leq e\leq|f|$, but his proof contains an error.  (See
the remarks after our lemma~\ref{lem-coprime-case}.)  To our
knowledge, the question of infiniteness of $G$ remains open for every
$G$ not treated in this paper, with two exceptions.   Jabara has informed the author that
he used the Knuth-Bendix algorithm in MAGNUS to find confluent
rewriting systems for $BS(2,3;2,3;2,3)$ and $BS(3,4;3,4;3,4)$, and
then counted the language of irreducible words to show the groups are
infinite.  

Not only does $G$ collapse in the situation of theorem~\ref{thm-non-devlopability}, but we
can say a great deal about what it collapses to.  And with no more work,
we also get information about the finite quotients of $G$ in
many cases not covered by theorem~\ref{thm-non-devlopability}.

\begin{theorem}
\label{thm-universal-quotient}
Suppose $(a,b)=(c,d)=(e,f)=1$ and none of the three pairs is
$(\pm1,\pm1)$.  Then there exists a quotient $Q=Q(a,b;c,d;e,f)$ of
$G(a,b;c,d;e,f)$ which is universal among all quotients in which $x$,
$y$ and $z$ have finite order; that is: any such quotient factors
through $Q$.  Furthermore, $Q$ is finite and solvable, with
its commutator subgroup $Q'$ nilpotent of class${}\leq2$.  Finally,
if any of $a,\dots,f$ is $1$ then $G=Q$.
\end{theorem}

This immediately implies Post's result \cite{Post} that
$G(a,a+1;b,b+1;1,2)$ is trivial, since it is a solvable group with
trivial abelianization.  In section~\ref{sec-finite-solvable-groups}
we provide more detailed information, like a formula for the
order of $Q$, exact up to a divisor of $(b-a)^2(d-c)^2(f-e)^2$, and
a result showing that $Q'$ is usually abelian, not just nilpotent.
But $Q'$ is not always abelian: a calculation using GAP \cite{GAP}
shows that $Q(1,4;1,4;1,4)'$ is nonabelian.

The special cases in theorem~\ref{thm-non-devlopability} indicate special behavior when $b=\pm a$, $d=\pm c$
or $f=\pm e$.  This reflects properties of the Baumslag-Solitar
groups
\begin{equation}
\label{eq-def-of-BS-group}
BS(a,b):=\bigl\langle x,y\bigm|(x^a)^y=x^b\bigr\rangle,
\end{equation}
which we recall here to help orient the reader.
First, $BS(1,\pm1)=\Z\semidirect\Z$, the quotient $\Z$ acting on the normal
subgroup $\Z$ trivially or by $\{\pm1\}$. Second, $BS(1,n\neq\pm1)$ is
$\Z[\frac{1}{n}]\semidirect\Z$, the a generator of the quotient $\Z$ acting on
$\Z[\frac{1}{n}]$ by multiplication by $n$.  Here $\Z[\frac{1}{n}]$
means the subring of $\Q$, or rather the underlying abelian group.
Third, if $(a,b)=1$ and $a,b\notin\{\pm1\}$ then $BS(a,b)$ contains
nonabelian free groups, and is non-Hopfian, non-residually-finite, and
non-linear \cite{BS}.  Finally, if $a$ and $b$ have a common divisor $l$
then $BS(a,b)$ is an amalgamated free
product of $BS(a/l,b/l)$ and $\Z$.

I am very grateful to E. Jabara for pointing me toward the older
literature on these groups, most of which I was unaware of.

\section{The relatively prime case}

In this section we will prove theorem~\ref{thm-non-devlopability} in the special case that
$(a,b)=(c,d)=(e,f)=1$.  This is the basis for the general proof in the
next section.  Our convention for conjugation is that $x^y=\ybar xy$,
where $\ybar$ means $y^{-1}$.  Also, since some superscripts get very
complicated, we sometimes write $x\power{y}$ for $x^y$.  

Our first step is to find the key relation that makes the triangles of
theorem~\ref{thm-non-devlopability} collapse; the exact form of the
relation is not so important---the key is that some power of $x$ lies
in $\langle y,z\rangle$.  The restriction to $a,\dots,f>0$ is minor,
as we will see in the proof of lemma~\ref{lem-x-has-finite-order}.

\begin{lemma}
\label{lem-killer-relation}
Suppose $0<a\leq b$, $0<c\leq d$ and $1=e\leq f$.  Then for any
$R,S,T>0$, the relation
\begin{multline}
\label{eq-killer-relation}
x\Bigpower{Tb^{Sc^R}\Bigl[b^{S(d^R-c^R)}-a^{S(d^R-c^R)}\Bigr]}\\
=\zbar\Bigpower{Rf^{Ta^{S(d^R-c^R)}b^{Sc^R}}}
\ybar^{Sc^R}
z\Bigpower{Rf^{Ta^{Sd^R}}}
y^{Sd^R}
\end{multline}
holds in $G(a,b;c,d;e,f)$.
\end{lemma}

\begin{proof}
We will evaluate
$(x^P)\power{(y^Q)^{z^R}}$ in two different ways, where $P,Q,R$ are
integers having whatever divisibility properties are needed for the
following calculation to make sense.  The underlines indicate where
changes occur.
\begin{align*}
(x^P)\power{y^{Q(d/c)^R}}
&{}=
\zbar^R\ybar^Q\, \underline{z^R x^P}\,\zbar^R y^Q z^R\\
x\power{P(b/a)^{Q(d/c)^R}}
&{}=
\zbar^R\,\underline{\ybar^Qx^P}\,z\power{R(f/e)^P}\underline{\zbar^R y^Q z^R}\\
&{}=
\underline{\zbar^Rx\power{P(b/a)^Q}}\,\ybar^Qz\power{R(f/e)^P}y^{Q(d/c)^R}\\
&{}=x\power{P(b/a)^Q}\zbar\Bigpower{R(f/e)^{P(b/a)^Q}}\\
&\kern80pt
\cdot\ybar^Qz\power{R(f/e)^P}y^{Q(d/c)^R}.\\
\end{align*}
(The second line uses $z^Rx^P=x^P(z^R)^{x^P}$ and the fourth line is
similar, while the third line uses $\ybar^Qx^P=(x^P)^{y^Q}\ybar^Q$.)
The restrictions on $P$, $Q$ and $R$ come from considerations like
this: $(y^Q)^{z^R}=y^{Q(d/c)^R}$ follows from $(y^c)^z=y^d$ provided
that $Q$ is divisible by $c^R/(c,d)^{R-1}$ if $R>0$, or by
$d^{|R|}/(c,d)^{|R|-1}$ if $R<0$.  
The full set of conditions for the calculation to make sense are $c^R|Q$,
$a^Q|P$ and $a^{Q(d/c)^R}{\bigm|}P$.  Since $Q(d/c)^R\geq Q$, the third condition
implies the second.  We obtain \eqref{eq-killer-relation} by taking $Q=Sc^R$ and
$P=Ta^{Sd^R}$. 

(One can show that there are no
solutions for $P,Q,R$ unless at least one of $b/a$, $d/c$ and $f/e$ is
an integer or reciprocal integer.  Each of these is${}\geq1$, so we
ignore the case of reciprocal integers.  Also, if only one of $b/a$,
$d/c$ and $f/e$ is an integer, then it must be $f/e$.  This would
explain a hypothesis $e|f$, and the stronger assumption $e=1$
because it is enough for our applications.)

\end{proof}

\begin{lemma}
\label{lem-x-has-finite-order}
Suppose $b\neq\pm a$ and $d\neq\pm c$.  Then $x$ has finite order in
$G(a,b;c,d;1,f)$.  If $f\neq\pm1$ then $y$ and $z$ also have finite order.
\end{lemma}

\begin{proof}
Note that $x^2$, $y^2$, $z^2$ satisfy the relations of
$G(a^2,b^2;c^2,d^2;e^2,f^2)$, so it suffices to treat the case
$a,\dots,f>0$.  So we may take $0<a<b$, $0<c<d$ and $1=e\leq f$
without loss.

First suppose $f>1$.  Take $R=S=T=1$ and write the relation \eqref{eq-killer-relation}
as $x^A=\zbar^B\ybar^C z^D y^E$.  Being a word in $y,z$, $x^A$ 
conjugates some power of $y$ to another power.  Namely,
$$
\bigl(y^{d^B}\bigr)^{x^A}
=
\bigl(y^{d^B}\bigr)^{\zbar^B\ybar^C z^D y^E}
=
\bigl(y^{c^B}\bigr)^{z^Dy^E}
=
y^{c^{B-D}d^D}.
$$
(The third equality is valid because $B>D$.)   We write this relation
as $(y^g)^X=y^h$ where $g=d^B$, $h=c^{B-D}d^D$ and $X=x^A$.  Now we
apply the relation $z^{f^A}=Z^X$ to a large power of $y$.  The
conjugate of $y\power{hc^{f^A}}$ by $z^{f^A}$ is $y\power{hd^{f^A}}$.
So we have
\begin{align*}
y\power{hd^{f^A}}
&{}=
\Xbar\zbar X\, y\power{hc^{f^A}}\,\Xbar z X\\
&{}=
\Xbar\zbar\, y\power{gc^{f^A}}\,zX\\
&{}=
\Xbar\, y\power{gc^{f^A-1}d}\, X\\
&{}=
y\power{hc^{f^A-1}d}.
\end{align*}
We conclude that $y$ has order dividing $hd(d^{f^A-1}-c^{f^A-1})$.
This is a nontrivial relation provided $d^{f^A-1}\neq c^{f^A-1}$.
Since $d>c>0$, the relation is nontrivial provided $f^A\neq1$.  Since
$f>1$, the relation is nontrivial provided $A\neq0$.  Recall that $A$
is the exponent on the left side of \eqref{eq-killer-relation} and
that $b>a>0$ and $d>c>0$, so $A\neq0$.  Therefore $y$ has finite order, say
$y^n=1$.  Now the relation $x^{a^n}=(x^{a^n})^{y^n}=x^{b^n}$ implies
that $x$ has finite order (since $b>a>0$), and repeating this argument
shows that $z$ also has finite order.

If $f=1$ then $[x,z]=1$ and the computation is similar but much
easier.  We conjugate $x^{a^d}$ by the relation $(y^c)^z=y^d$, which
leads to $x\power{a^{d-c}b^c}=x\power{b^d}$.  Since $b>a>0$ and
$d>c>0$, this is a nontrivial relation, so $x$ has finite order.
(Remark: $y$ and $z$ have infinite order, since adjoining the
relation $x=1$ reduces $G$ to $BS_{yz}(c,d)$.)
\end{proof}

\begin{lemma}
\label{lem-coprime-case}
Assume the hypotheses of Theorem~\ref{thm-non-devlopability} and that
$(a,b)=(c,d)=(e,f)=1$.  Then the triangle of groups there is
developable in and only in the following cases:
\begin{align}
\label{eq-coprime-special-case-1}
&G(1,-1;1,-1;1,-1)\\
\label{eq-coprime-special-case-2}
&G(a,b;1,1;1,1)\\
\label{eq-coprime-special-case-3}
&G(a,b;1,1;1,-1)\hbox{, $a,b$ odd.}
\end{align}
\end{lemma}

\begin{proof}
The hypothesis that
one of $a,\dots,f$ divides its partner says that one is $\pm1$, say
$e=1$ without loss.  If $b\neq\pm a$ and $d\neq\pm c$ then
lemma~\ref{lem-x-has-finite-order} shows that $G$ is not developable.  So suppose $b=\pm a$
or $d=\pm c$.  Because of the relative primality, we are in one of the
cases $G(1,\pm1;1,\pm1;1,f)$,  $G(1,\pm1;c,d\neq\pm c;1,f)$ and
$G(a,\discretionary{}{}{}b\neq\pm a;\discretionary{}{}{}1,\discretionary{}{}{}\pm1;\discretionary{}{}{}1,\discretionary{}{}{}f)$.  In the last two cases, when
$f\neq\pm1$, $G$ is non-developable by lemma~\ref{lem-x-has-finite-order} (after cyclically
permuting the variables).  All remaining cases are now special cases
of $G(a,b;1,\pm1;1,\pm1)$, after cyclic permutation of the variables.
To begin with, $G(a,b;1,1;1,1)=BS_{xy}(a,b)\times\Z_z$ is obviously
developable. 

Next, in $G(a,b;1,1;1,-1)$, $\langle z\rangle$ is normal, with
quotient $BS(a,b)$, so $\langle x,y\rangle$ is a complement.  So $G=\langle z\rangle\semidirect
BS_{xy}(a,b)$, with $y$ fixing $z$ and $x$ inverting it.  If $a$ and
$b$ have different parities then $(x^a)^y=x^b$ forces $z^2=1$, so $G$
is not developable.  On the other hand, if $a$ and $b$ have the same
parity (so both are odd, because $(a,b)=1$), then $\Z_z\semidirect
BS_{xy}(a,b)$ satisfies all the relations of $G$, hence equals it.  So
$G(a,b;1,1;1,-1)$ is developable if and only if $a$ and $b$ are both
odd.

Next, in $G(a,b;1,-1;1,1)$, $\langle x,y\rangle$ is normal, with a
complementary $\Z$ generated by $z$.  Conjugating $(x^a)^y=x^b$ by $z$
gives $(x^a)^{\ybar}=x^b$, so $x^{a^2}=(x^{a^2})^{y\ybar}=x^{b^2}$.
Therefore $x$ has finite order unless $b=\pm a$, which by $(a,b)=1$
leaves us with $G(1,\pm1;1,-1;1,1)$, which we treated in the previous
case. 

Finally we consider $G(a,b;1,-1;1,-1)$.  Because $x^2,y,z$ satisfy the
relations of $G(a,b;1,-1;1,1)$, the previous case shows that only
$G(1,\discretionary{}{}{}\pm1;\discretionary{}{}{}1,\discretionary{}{}{}-1;\discretionary{}{}{}1,\discretionary{}{}{}-1)$
can be developable.  The $+1$ case has already been treated, leaving
only $G(1,-1;1,-1;1,-1)$, whose developability is due to Neumann
\cite[\S5]{Neumann}.  Observe that $x^2$, $y^2$ and $z^2$ generate a
normal abelian group $A$, with quotient $(\Z/2)^3$.  To see that $G$
is developable, it suffices to prove $A\iso\Z^3$.  This can be done by
representing $G$ by isometries of $\R^3$, with $x$ acting by
$(X,Y,Z)\mapsto(X+1,Y,-Z)$ and the other generators' actions defined
similarly.  In fact this action on $\R^3$ is free, realizing $G$ as
the fundamental group of a Euclidean 3-manifold.
\end{proof}

\begin{remark}
The group $G=G(2,3;2,3;2,4)$ is not developable, because $x,y,z^2$
satisfy the relations of $G(2,3;4,9;1,2)$, and the latter is
non-developable by lemma~\ref{lem-coprime-case}.  This group is a counterexample to
the main result (theorem~4.4) of \cite{Neumann-broken}.  Neumann's argument relies
on a complicated inductive definition of an action of $G$ on a set of
``normal matrices''.  Unfortunately, his operator $\rho(b^{-1})$
doesn't preserve the set of normal matrices: the right hand side of
(3.53) is never a normal matrix because it violates (2.35) or (2.36),
depending on the sign of $\gamma(n)$.  (His proof of the nonexistence of
finite quotients of $G(a,a+1;c,c+1;e,e+1)$ is correct.)
\end{remark}

\section{The General Case}
\label{sec-general}

In this section we derive theorem~\ref{thm-non-devlopability} from the coprime case
established in lemma~\ref{lem-coprime-case}.  The key idea is the following; consider
$G:=G(a,b;c,d;e,f)$ and suppose $l>0$ is a common divisor of $a$ and $b$.
Then the elements $X=x^l$, $y$ and $z$ satisfy the relations of
$H:=G(a/l,b/l;c,d;e^l,f^l)$.  Because of this change of variables, we
will sometimes refer to \eqref{eq-deg-of-G} as $G_{xyz}(a,b;c,d;e,f)$
and \eqref{eq-def-of-BS-group} as $BS_{xy}(a,b)$.  In this notation,
$G_{xyz}(a,b;c,d;e,f)$ is the direct limit of the diagram
\begin{equation}
\label{eq-pushout-in-nonprime-case}
G_{Xyz}(a/l,b/l;c,d;e^l,f^l)\leftarrow
BS_{zX}(e^l,f^l)
\to
BS_{zx}(e,f).
\end{equation}
We also sometimes write $\Z_x$ for a copy of $\Z$ with generator $x$.
The right homomorphism of \eqref{eq-pushout-in-nonprime-case} is always
injective; to see this, one may use the standard form for words in an
HNN extension.  In good cases, the left homomorphism is also
injective, so that $G$ is an amalgamated free product of $H$ and $BS_{zx}(e,f)$.  When this
holds, we may reasonably hope to relate the developability of $G$ to
that of $H$.  This hope is realized in the
following lemma.

\begin{lemma} In the notation just established,
\label{lem-developability-criteria}
\begin{enumerate}
\item
\label{item-G-developability-implies-H-developability}
If $G$ is developable then so is $H$.
\item
\label{item-H-developability-almost-implies-G-developability}
Suppose that $H$ is developable and that $\langle
X,y\rangle\cap\langle z,X\rangle=\langle X\rangle$.  Then $G$ is developable.
\end{enumerate}
\end{lemma}

\begin{proof}
\eqref{item-G-developability-implies-H-developability}  If $H$ is not
developable then $BS_{X,y}(a/l,b/l)$, $BS_{yz}(c,d)$ or
$BS_{zX}(e^l,f^l)$ fails to inject into $H$.  Since these are
subgroups of $BS_{xy}(a,\discretionary{}{}{}b)$, $BS_{yz}(c,d)$ and $BS_{zx}(e,f)$, at
least one of these latter three fails to inject into $G$.

\eqref{item-H-developability-almost-implies-G-developability} The left arrow of
\eqref{eq-pushout-in-nonprime-case} is injective, by the definition of developability of $H$.  So
\eqref{eq-pushout-in-nonprime-case} expresses $G$ as a free product with amalgamation.  Since
$BS_{zx}(e,f)$ is a factor in this product, it injects into $G$.
Also, $BS_{yz}(c,d)$ injects into $H$ by developability, and then
injects into $G$ since $H$ does.  So it remains to check the
injectivity of $BS_{xy}(a,b)$ into $G$.

We use the following assertion, whose proof is an easy
exercise using the standard form for words in an amalgamated free
product.  Suppose we are given a commutative diagram of inclusions of
groups 
$$
\begin{CD}
A@<<<I@>>>B\\
@AAA @AAA @AAA\\
C@<<< J @>>> D;
\end{CD}
$$
then $I\cap C=J=I\cap D$ implies that the natural map $C*_JD\to A*_IB$
is injective.
The hypothesis in \eqref{item-H-developability-almost-implies-G-developability} is exactly what is needed to apply this to
the diagram
$$
\begin{CD}
H @<<< BS_{zX}(e^l,f^l) @>>> BS_{zx}(e,f)\\
@AAA @AAA @AAA\\
BS_{Xy}(a/l,b/l)@ <<< \Z_X @>>> \Z_x.
\end{CD}
$$ The amalgamation of the bottom row is $BS_{xy}(a,b)$ and that of
the top is $G$.  So the former injects into the latter and the proof
is complete.
\end{proof}

In order to deduce the developability of $G$ from that of $H$, we must
verify the condition in \eqref{item-H-developability-almost-implies-G-developability}.  We will prove this in
lemma~\ref{lem-developable-implies-vertex-groups-intersect-well}, by an
argument that requires understanding certain centralizers in $H$:

\begin{lemma}
\label{lem-centralizers-in-BS}
In $BS_{xy}(a,b)$, the centralizer of $y^n$ is 
\begin{enumerate}
\item
\label{item-centralizer-Z-by-Z}
$\langle x^a,y\rangle$ if $a=b$, or if $a=-b$ and $n$ is even;
\item
\label{item-centralizer-Z}
$\langle y\rangle$ otherwise.
\end{enumerate}
\end{lemma}

\begin{proof}
This is an exercise using the standard form for words in an HNN
extension.  Or one can apply the last part of the theorem stated on
pp. 350--351 of \cite{KS}.
\end{proof}

Now we verify the condition in
lemma~\ref{lem-developability-criteria}\eqref{item-H-developability-almost-implies-G-developability}.
Part \eqref{item-required-for-simultaneous-induction} of the following
lemma is needed for the inductive argument, but nowhere else.  The
important conclusion is
\eqref{item-vertex-groups-intersect-as-expected}.

\begin{lemma}
\label{lem-developable-implies-vertex-groups-intersect-well}
Suppose $G$ is developable.  Then
\begin{enumerate}
\item
\label{item-vertex-groups-intersect-as-expected}
$\langle{x,y}\rangle\cap\langle y,z\rangle=\langle y\rangle$ and
similarly for cyclic permutations of $x,y,z$;
\item
\label{item-required-for-simultaneous-induction}
if $|a|=|b|$, $|c|=|d|$ and $|e|=|f|$ then some powers of $x$, $y$ and
$z$ generate a group $\Z^3$.
\end{enumerate}
\end{lemma}

\begin{proof}
Suppose $G$ were a counterexample, with $|a|+\cdots+|f|$ minimal.
If it is conclusion \eqref{item-required-for-simultaneous-induction}
that fails for $G$, then $b=\pm a$, $d=\pm c$ and $f=\pm e$.  We
cannot have $a,\dots,f\in\{\pm1\}$, because then we would be in one of
the special cases $G=G(1,\pm1;1,\pm1;1,\pm1)$, for which the lemma can
be checked directly.  (The only interesting case is
$G(1,-1;1,-1;1,-1)$, for which see the proof of lemma~\ref{lem-coprime-case}.)  So
suppose $a>1$, so that $G$ is the pushout of the diagram
\begin{equation}
\label{eq-foo}
G_{Xyz}(1,\pm1;c,d;e^a,f^a)\leftarrow BS_{zX}(e^a,f^a)\rightarrow
BS_{zx}(e,f). 
\end{equation}
The developability of $G$ implies that of the left term $H$
(lemma~\ref{lem-developability-criteria}\eqref{item-G-developability-implies-H-developability}),
so \eqref{eq-foo} expresses $G$ as a free product with amalgamation,
so $H$ injects into $G$.  Now applying the inductive hypothesis to
$H$, we see that some powers of $X,y,z$ generate a group $\Z^3$.
Since $X$ is a power of $x$, we have proven
\eqref{item-required-for-simultaneous-induction}.

So it must be \eqref{item-vertex-groups-intersect-as-expected} that
fails.  Then $\langle{x,y}\rangle\cap\langle y,z\rangle$ is strictly
larger than $\langle y\rangle$, so take $w$ to be an element in the
intersection but not in $\langle y\rangle$.  Since $w\in\langle
y,z\rangle$, it conjugates some power of $y$ to another power
(possibly the same), say $(y^m)^w=y^n$.  On the other hand, since
$w\in\langle x,y\rangle$, we see that $y^m$ and $y^n$ are conjugate in
$\langle x,y\rangle= BS_{xy}(a,b)$.  This forces $m=n$, so that $w$
centralizes some power of $y$.  Since $w\notin\langle y\rangle$,
lemma~\ref{lem-centralizers-in-BS} forces $a=\pm b$ and
$w\in\langle{x^a,y}\rangle=\langle x^a\rangle\semidirect\langle
y\rangle$.  Any subgroup of this $\Z\semidirect\Z$ that strictly
contains $\langle y\rangle$ must contain a power of $x$.  Therefore
$\langle{x,y}\rangle\cap\langle y,z\rangle$ contains a power of $x$;
we may even suppose without loss of generality that $w$ is a power of
$x$.

As a power of $x$, $w$ conjugates some power of $z$ to
another, say $(z^p)^w=z^q$.  We now essentially repeat the argument just
used:  since $w\in\langle y,z\rangle=BS_{yz}(c,d)$, we must have
$p=q$, and this forces $f=\pm e$.  Also, since $w$ centralizes a power
of $z$ and is not in $\langle y\rangle$, the centralizer of $z$ in
$BS_{yz}(c,d)$ must be larger than $\langle y\rangle$, which forces
$c=\pm d$ by lemma~\ref{lem-centralizers-in-BS}.

We have proven that $a=\pm b$, $c=\pm d$, $e=\pm f$ and that some
power of $x$ lies in the centralizer of a power of $z$ in
$BS_{yz}(c,\pm c)$, which has structure $\langle
y^c\rangle\semidirect\langle z\rangle$.  
But this contradicts the fact that some powers of $x,y,z$ generate a copy
of $\Z^3$, by \eqref{item-required-for-simultaneous-induction}. 
\end{proof}

We summarize our results so far as:

\begin{lemma}
\label{lem-G-developable-iff-H-is}
$G$ is developable if and only if $H$ is.
\end{lemma}

\begin{proof}
We have already shown that developability of $G$ implies that of $H$.
For the converse, we apply
lemma~\ref{lem-developable-implies-vertex-groups-intersect-well} to
$H$, and then conclusion
\eqref{item-vertex-groups-intersect-as-expected} of that lemma allows
us to apply lemma~\ref{lem-developability-criteria} and deduce $G$'s
developability.
\end{proof}

\begin{corollary}
\label{cor-nonprime-developability}
Write 
$$
(a,b;c,d;e,f)=(Al,Bl;Cm,Dm;En,Fn),
$$
where $l,m,n>0$ and
$(A,B)=(C,D)=(E,F)=1$.  Then $G(a,\discretionary{}{}{}b;\discretionary{}{}{}c,\discretionary{}{}{}d;\discretionary{}{}{}e,\discretionary{}{}{}f)$ is developable if and
only if $G(A^m,B^m;C^{n^l},D^{n^l};E^l,F^l)$ is.
\end{corollary}

\begin{proof}
Consider the following four groups:
\begin{gather}
G(Al,Bl;Cm,Dm;En,Fn)\\
G(A,B;Cm,Dm;(En)^l,(Fn)^l)\\
G(A^m,B^m;C,D;(En)^l,(Fn)^l)\\
G(A^m,B^m;C^{n^l},D^{n^l};E^l,F^l).
\end{gather}
By lemma~\ref{lem-G-developable-iff-H-is}, each is developable if and only if the previous
one is.
\end{proof}

\begin{proof}[Proof of theorem~\ref{thm-non-devlopability}:]
We suppose without loss that $a,c,e>0$.  If $b=\pm a$, $d=\pm c$ and
$f=\pm e$ then we are in case
\begin{equation}
\label{eq-+-+-+-}
G=G(a,\pm a;c,\pm c;e,\pm e)
\end{equation}
and corollary~\ref{cor-nonprime-developability} and lemma~\ref{lem-coprime-case} imply that $G$ is developable.
If two of the equalities $b=\pm a$, $d=\pm c$, $f=\pm e$ fail, then
the corollary and lemma prove $G$ non-developable.  The remaining case
is when exactly one of the equalities fails, so suppose $b\neq a$,
$d=\pm c$, $f=\pm e$.  We take $l,m,n,A,\dots,F$ as in
corollary~\ref{cor-nonprime-developability}.  Since $a,c,e>0$ we have
$A,C,E>0$.  By that corollary, $G$ is developable if and only if
$G(A^m,B^m;C^{n^l},D^{n^l};E^l,F^l)$ is, which can be determined using
the relatively-prime case, lemma~\ref{lem-coprime-case}.  
So developability is equivalent to
$(A^m,B^m;C^{n^l},D^{n^l},E^l,F^l)$ being equal to
\begin{align}
\label{eq-case-1}
&(A^m,B^m;1,1;1,-1)\hbox{ with $A^m$ and $B^m$ odd}\\
\label{eq-case-2}
\hbox{or }&(A^m,B^m;1,1;1,\-1).
\end{align}
In either case, we know $C=E=1$ because $C,E>0$.

In case \eqref{eq-case-1}, $F^l=-1$ is equivalent to $F=-1$ and $l$
odd, and of course the oddness of $A^m$ and $B^m$ is equivalent to the
oddness of $A$ and $B$.  The condition $D^{n^l}=1$ is equivalent to:
either $D=1$, or else $D=-1$ and $n$ is even.  So we have
\begin{align*}
(A,B;C,D;E,F)
={}&(A,B;1,\-1;1,-1)\hbox{ with $A$, $B$, $l$
    odd},\\
\hbox{or }&(A,B;1,-1;1,-1)\hbox{ with $A,B,l$ odd and $n$ even;}
\end{align*}
note that $n=e$.  This is equivalent to 
\begin{align}
\label{eq-ab-odd+-}
G={}&G(a,b;c,\-c;e,-e)\hbox{, with $a,b$ odd}\\
\label{eq-ab-odd--}
\hbox{or }&G(a,b;c,-c;e,-e)\hbox{, with $a,b$ odd and $e$ even}.
\end{align}

In case \eqref{eq-case-2}, the treatment of $D^{n^l}=1$ is as before,
and $F^l=1$ is equivalent to: either $F=1$, or else $F=-1$ and $l$ is
even.  So we have
\begin{align*}
(A,B;C,D;E,F)
={}&(A,B;1,\-1;1,\-1),\\
\hbox{or }&(A,B;1,\-1;1,-1)\hbox{ with $l$ even,}\\
\hbox{or }&(A,B;1,-1;1,\-1)\hbox{ with $n$ even,}\\
\hbox{or }&(A,B;1,-1;1,-1)\hbox{ with $l$ and $n$ even;}
\end{align*}
again $n=e$.  This is equivalent to
\begin{align}
\label{eq-ab}
G={}&G(a,b;c,\-c;e,\-e)\\
\label{eq-ab+-even}
\hbox{or }&G(a,b;c,\-c;e,-e)\hbox{, with $a,b$ even}\\
\label{eq-ab-+e-even}
\hbox{or }&G(a,b;c,-c;e,\-e)\hbox{, with $e$ even}\\
\label{eq-ab--even}
\hbox{or }&G(a,b;c,-c;e,-e)\hbox{, with $a,b,e$ even}.
\end{align}

Now, \eqref{eq-ab-odd+-} and \eqref{eq-ab+-even} together correspond
to \eqref{eq-special-case-3} in the statement of the theorem, and
\eqref{eq-ab-odd--} and \eqref{eq-ab--even} correspond to
\eqref{eq-special-case-5}.  Also, \eqref{eq-ab} and
\eqref{eq-ab-+e-even} correspond to \eqref{eq-special-case-2} and
\eqref{eq-special-case-4}, and
\eqref{eq-special-case-2}--\eqref{eq-special-case-5} contain every
case of \eqref{eq-+-+-+-} except $G(a,-a;c,-e;e,-e)$, which we listed
as \eqref{eq-special-case-1}.
\end{proof}

\section{Finite Solvable Groups}
\label{sec-finite-solvable-groups}

We have shown that $G=G(a,b;c,d;e,f)$ is non-developable under fairly
mild conditions, and in this section we study just how much $G$
collapses.  We first prove
theorem~\ref{thm-universal-quotient}, which often says that
$G$ is a finite solvable group.  We assume the hypotheses of
theorem~\ref{thm-universal-quotient} throughout this section, and
without loss we suppose $a<b$, $c<d$, $e<f$.  It is convenient to
define $X=x^{b-a}$, $Y=y^{d-c}$ and $Z=z^{f-e}$.

\begin{lemma}
\label{lem-existence-of-universal-quotient}
The relation
\begin{equation}
\label{eq-order-of-x-in-finite-quotient}
x\power{(b-a)^2(b^{d-c}-a^{d-c})}=1
\end{equation}
and its cyclic permutations hold in any quotient of $G$ in which $x$,
$y$ and $z$ have finite order.  In particular, $G$ has a universal
quotient $Q$ in which $x$, $y$ and $z$ have finite order, which
is got by imposing these relations.
\end{lemma}

\begin{proof}
Suppose $\Gbar$ is a quotient of $G$ in which $x,y,z$ have finite
order, and write $n$ for the order of $x$.   The orders of $x^a$ and
$x^b$ are $n/(n,a)$ and $n/(n,b)$, which are equal since $x^a$
and $x^b$ are conjugate.  Since $(a,b)=1$, this forces
$(n,a)=(n,b)=1$.  Therefore $\langle x^a\rangle=\langle
x^b\rangle=\langle x\rangle$, so $y$ normalizes $\langle x\rangle$.
Similarly, $z$ normalizes $\langle y\rangle$ and $x$ normalizes
$\langle z\rangle$.  

Now let $H$ be the subgroup generated by all the
$y^{x^i}$, $i\in\Z$.  We have $H=\langle y,x^{b-a}\rangle$, since
$y^{x^a}=yx^{a-b}$ and $x^a$ generates $\langle x\rangle$.  Obviously
$x$ and $y$ normalize $H$.  And the fact that $z$
normalizes $\langle y\rangle$ implies that $z^{x^i}$ normalizes
$\langle y^{x^i}\rangle$.  Since $\langle z\rangle=\langle
z^{x^i}\rangle$ normalizes every $\langle y^{x^i}\rangle$, it
normalizes $H$.  So $H$ is normal in $\Gbar$.

Next, the commutator subgroup $H'$ is $\langle x^{(b-a)^2}\rangle$,
which is characteristic in $H$, hence normal in $\Gbar$.  Now, the
automorphism group of a cyclic group is abelian, so every commutator
acts trivially, in particular $y^{d-c}$.  This implies
\eqref{eq-order-of-x-in-finite-quotient} and similarly for $y$ and
$z$.
\end{proof}

\begin{lemma}
\label{lem-conjugation-relations}
Let $\a$ be a solution of $\alpha a=1$ modulo
$(b-a)^2(b^{d-c}-a^{d-c})$.  Then
\begin{align}
\label{eq-y-conjugated-by-x}
y^x&{}= y x\power{-\a(b-a)}=yX^{-\a}\\
\label{eq-Y-conjugated-by-x}
Y^x&{}= Y X\power{-\a^{d-c}\frac{b^{d-c}-a^{d-c}}{b-a}}\\
\label{eq-Y-conjugated-by-X}
Y^X&{}= Y X\power{-\a^{d-c}(b^{d-c}-a^{d-c})}.
\end{align}
\end{lemma}

\begin{proof}
The key property of $\a$ is that $(x^a)^\a=1$.  We may rewrite
$(x^a)^y=x^b$ as $y^{x^a}=yx^{a-b}$.  Conjugating $y$ by $x^a$, $\a$
many times, gives \eqref{eq-y-conjugated-by-x}.  For
\eqref{eq-Y-conjugated-by-x} we compute
\begin{align*}
Y^{x^a}
&{}=
(y^x)^{d-c}=\bigl(y\Xbar^\a\bigr)^{d-c}\\
&{}=
Y\bigl(\Xbar^\a\bigr)^{y^{d-c-1}}\bigl(\Xbar^\a\bigr)^{y^{d-c-2}}\cdots\bigl(\Xbar^\a\bigr)^{y^{0}}\\
&{}=
Y\bigl(\Xbar^{\a\a^{d-c-1}a^{d-c-1}}\bigr)^{y^{d-c-1}}\cdots\bigl(\Xbar^{\a\a^{d-c-1}a^{d-c-1}}\bigr)^{y^0}\\
&{}=
Y\Xbar\Bigpower{\a^{d-c}\Bigl(b^{d-c-1}+b^{d-c-2}a+\cdots+a^{d-c-1}\Bigr)}\\
&{}=
YX\Bigpower{-\a^{d-c}\,\frac{b^{d-c}-a^{d-c}}{b-a}}.
\end{align*}
Then \eqref{eq-Y-conjugated-by-X} follows by applying \eqref{eq-Y-conjugated-by-x} $b-a$ times.
\end{proof}

\begin{proof}[Proof of theorem~\ref{thm-universal-quotient}:]
We must show that $Q'$ is nilpotent of class${}\leq2$.  It follows
from \eqref{eq-Y-conjugated-by-x} and its cyclic permutations that
$\langle X,Y,Z\rangle$ is normal in $Q$.  Since adjoining the
relations $X=Y=Z=1$ abelianizes $Q$, we see that $\langle
X,Y,Z\rangle=Q'$.  Then \eqref{eq-Y-conjugated-by-X} shows that $[X,Y]$ lies in $\langle
X^{b-a}\rangle$.  We saw in the proof of lemma~\ref{lem-existence-of-universal-quotient} that $\langle
X^{b-a}\rangle$ is central in $Q'$.  Together with the cyclic
permutations of this argument, we have proven that $[Q',Q']$ is
central in $Q'$, as desired.

For the final assertion of the theorem, just use
lemma~\ref{lem-x-has-finite-order}, which assures us that $x,y,z$ have
finite order in $G$, so $G$ must equal $Q$.
\end{proof}

Jabara \cite{Jabara} proved the stronger result that $Q''$ is central in
$Q$, not just in $Q'$.  He treated only the case $a=c=e=1$, but there
is no loss of generality because $\langle x\rangle=\langle
x^a\rangle=\langle x^b\rangle$ in $Q$, and similarly for $y$ and $z$.  

$Q'$ is abelian in almost all cases.  The easiest way to address this
question is to work one prime at a time, since the nilpotence of $Q'$
implies that $Q'$ is the direct product of its Sylow subgroups.  So
for a prime $p$ we define $Q_p$ as the quotient of $Q$ by all the
Sylow subgroups of $Q'$ except for the one associated to $p$.
Obviously, $Q'$ is abelian if and only if every $Q_p':=(Q_p)'$ is.

We said in the introduction that $Q(1,4;1,4;1,4)'$ is nonabelian.  We
found this using GAP \cite{GAP}, but simply entering the presentation
led to memory overflow during coset enumeration.  Adjoining the
relations $x^{81}=y^{81}=z^{81}=1$, which reduce $G$ to $Q_3$, let GAP
perform the computation almost instantly.

\begin{lemma}
\label{lem-p-doesnt-divide-all-three}
Unless $p$  divides  $b-a$, $d-c$ and $f-e$,  $Q_p'$ is
abelian. 
\end{lemma}

\begin{proof}
Since a nonabelian $p$-group has noncyclic Frattini quotient, it
suffices to show that $Q_p/\Phi(Q_p)$ is cyclic.  This is an abelian
group with generators $X,Y,Z$ satisfying relations including
$pX=pY=pZ=0$ and $(b-a)X=(d-c)Y=(f-e)Z=0$, in addition notation.
Suppose $p\notdiv d-c$, so $Y=0$.  If $p\notdiv b-a$ then $X=0$ and
$Q_p/\Phi(Q_p)$ is generated by $Z$, hence cyclic.  So suppose
$p|b-a$.  Conjugating the relation $Y=0$ by $x$ yields
$$
Y-\a^{d-c}\frac{b^{d-c}-a^{d-c}}{b-a}X=0
\hbox{, hence }
\frac{b^{d-c}-a^{d-c}}{b-a}X=0.
$$
The hypotheses $p|b-a$ and $p\notdiv d-c$ imply that the $p$-part of
the numerator is the same as that of the denominator.  So this
relation implies $X=0$, and again $Q_p/\Phi(Q_p)$ is cyclic.
\end{proof}

\begin{corollary}
\label{cor-no-common-divisor}
If $b-a$, $d-c$ and $f-e$ have no common divisor then $Q'$ is abelian.
\qed
\end{corollary}







Mennicke \cite{Mennicke} gave an order formula for $G(1,t;1,t;1,t)$, and
Johnson and Robertson \cite{JR} gave an upper bound for the order of
$G(1,b;1,d;1,f)$.  In \cite{AA}, Albar and Al-Shuaibi improve this bound
and give a correction to Mennicke's paper.  It seems that the exact
order and structure of $Q_p'$ depend sensitively on the number times
$p$ divides $b-a$, $d-c$ and $f-e$.  We offer upper and lower
bounds on $|Q|$ that are fairly close to each other:

\begin{theorem}
\label{thm-order-of-Q}
Suppose $a<b$, $c<d$ and $e<f$.  Then the order of $Q$ is
\begin{multline*}
\bigl(b^{d-c}-a^{d-c}\bigr)
\bigl(d^{f-e}-c^{f-e}\bigr)
\bigl(f^{b-a}-e^{b-a}\bigr)\\
\times
\hbox{\rm a divisor of}\ 
(b-a)^2(d-c)^2(f-e)^2.
\end{multline*}
\end{theorem}

\begin{proof}
Killing $x$ reduces $Q$ to a group in which $y$ has order
$d^{f-e}-c^{f-e}$.  This shows that
$d^{f-e}-c^{f-e}$ divides 
$[\langle x,y\rangle:\langle x\rangle]$, hence
$[\langle X,y\rangle:\langle X\rangle]$.
Similarly, killing $z$ shows that the order of $x$ is divisible by
$b^{d-c}-a^{d-c}$, so the order of $X$ is divisible by
$\bigl(b^{d-c}-a^{d-c}\bigr)/(b-a)$.  And killing $y$ leaves a
group of order $\bigl(f^{b-a}-e^{b-a}\bigr)(b-a)$.  Putting all this
together shows that
$$
|Q|=[Q:\langle y,X\rangle]\cdot[\langle y,X\rangle:\langle X\rangle]\cdot
[\langle X\rangle:1]
$$
is divisible by $\bigl(b^{d-c}-a^{d-c}\bigr)
\bigl(d^{f-e}-c^{f-e}\bigr)
\bigl(f^{b-a}-e^{b-a}\bigr)$.

On the other hand, the structure of $Q$ as a polycyclic group shows
that $|Q|$ divides the product of the orders of $x$, $y$ and $z$.
Referring to \eqref{eq-order-of-x-in-finite-quotient} shows that $|Q|$ divides
$$
\bigl(b^{d-c}-a^{d-c}\bigr)
\bigl(d^{f-e}-c^{f-e}\bigr)
\bigl(f^{b-a}-e^{b-a}\bigr)
(b-a)^2(d-c)^2(f-e)^2.
$$
\end{proof}

\end{document}